\documentclass{article}
\usepackage{indentfirst, amsfonts, amsmath, amsthm,amssymb,amscd}
\usepackage[dvips]{graphicx}

\begin{document}
\title {\bf The scalar curvature and the biorthogonal curvature: A pinching problem}
\author { Ezio de Araujo Costa}
\date{}

\maketitle

\begin{abstract} { The famous pinching problem says that on a compact simply connected $n$-manifold if its sectional curvature satisfies  $K_{min} > \frac{1}{4}K_{max} > 0$, then the manifold is homeomorphic to the sphere. In [8, problem 12], S. T. Yau  proposed the following problem: If we replace $K_{max}$ by the scalar curvature, can we deduce similar pinching theorems? In our present note we
give an answer to this question in dimension $n = 4.$}

\end{abstract}
\noindent{\bf Mathematic subject classifications (2000): 53C25, 53C24. \noindent
}
\\
\\
Let $M = M^n$ be a connected differentiable $n$-manifold with  scalar curvature $s$ and let $K$ be the sectional curvature. For each $x\in M$ lets consider $P_1, P_2$ two mutually orthogonal 2-planes in the tangent space $T_xM$. The {\it biorthogonal (sectional) curvature } (see [2]) relative to $P_1$ and $P_2$  is the average $$K^\perp (
P_1, P_2) = \frac{K(P_1) + K(P_2)}{2}. \eqno [1.1]$$
If $n = 4$ and $P$  is a 2-plane in $T_xM$ we will write $$K^\perp (P) = \frac{K(P) + K(P^\perp) }{2}, \eqno [1.2]$$ where $P^\perp$ is the orthogonal complement of  $P$ in $T_xM$. The sum of two sectional curvatures on two orthogonal planes  appears in the works of W. Seeman [6] and M. H. Noronha [5]. In our presents article we consider  manifolds $M$ of dimension four and the following functions:

 $$K_1^\perp(x) = \textmd{min} \{K^\perp(P); P  \textmd{ is a 2- plain in } T_xM \}, \eqno [1.3]$$

  $$
  K_3^\perp(x) = \textmd{max}\{K^\perp (P); P \textmd{ is a 2- plain in } T_xM \}, \eqno [1.4]$$

 $$K_2^\perp (x)= \frac{s(x)}{4} - K_1^\perp(x) - K_3^\perp(x). \eqno [1.5]$$

\newpage

\begin{center} {\bf The biorthogonal curvature and the Weyl tensor }
\end{center}
Let $(M, g)$ be an oriented Riemannian 4-manifold.
For each $x\in M$ the bundle of two-forms $\Lambda^2 $ of $M$ splits
$\Lambda^2  = \Lambda^+ \bigoplus\Lambda^-$ into $\pm$-eigenspaces of the Hodge
$\star$-operator: $\Lambda^\pm = \{\alpha \in \lambda^2 ; \star\alpha = \pm\alpha \}$. The Weyl tensor $W$ is an endomorphism of
$\Lambda^2 $ such that $W = W^+\bigoplus W^-$, where $W^\pm : \Lambda^\pm
\longrightarrow \Lambda^\pm$ are self-adjoint with free traces and are called of the self-dual and anti-self-dual parts of $W$, respectively.
 Let $w_1^\pm \leq w_2^\pm \leq w_3^\pm$ be the eigenvalues of the tensors $W^\pm$, respectively. As was proved in [2],
$$K_1^\perp - \frac{s}{12} = \frac{w_1^+ + w_1^-}{2}, \eqno [1.6]$$
  $$K_2^\perp - \frac{s}{12} = \frac{w_2^+ + w_2^-}{2} \eqno [1.7]$$ and
 $$K_3^\perp - \frac{s}{12} = \frac{w_3^+ + w_3^-}{2} \eqno [1.8]$$

Based on a proposed question by Yau [8, problem 12], the authores of the article [3, page 16] proposed the following conjecture
\\
\\
{\it Let $(M, g)$ be a compact simply connected Riemannian n-manifold scalar curvature $s > 0$ and sectional curvature $K$.
If $K > \frac{s}{n(n+2)}$ on $M$ then $M$ is diffeomorphic to the sphere $\mathbb{S}^4$.}
 \\
 \\
  In dimension $n = 4$ we obtained the following
\\
\\
{\bf Theorem 1} - {\it Let $(M,g)$ be a compact oriented 4-manifold with scalar curvature $s > 0$. Let $K_1^\perp$ and $K_3^\perp$ be the biorthogonal curvatures given by [1.3] and [1.4], respectively.
If $K_1^\perp \geq \frac{s}{24}$ on $M$ or $K_3^\perp \leq \frac{s}{6}$ on $M$ then we have
 \\
 \\
 (1) $M$ is diffeomorphic to a connected sum $\mathbb{S}^4 \sharp(\mathbb{R} \times \mathbb{S}^3)/G_1\sharp ......\sharp(\mathbb{R} \times \mathbb{S}^3)/G_n$, where the $G_i$ are discrete subgroup of the isometry group of $\mathbb{R} \times \mathbb{S}^3)$ or
\\
\\
(2) $(M, g)$ is conformal to a complex projective space $\mathbb{CP}^2$ with the Fubini-Study metric or
\\
\\
(3) The  universal covering of $M$ is isometric to product $\mathbb{R} \times N^3$, where $N^3$ is diffeomorphic  to $\mathbb{S}^3$. }
\\
\\

{\bf Remark 1.1}- Compare the above Theorem 1 to Theorem 1.1 in [3] and the Conjecture A in page 17 of [3].
\newpage

\begin{center} {\bf Proof of Theorem 1}
\end{center}
Let $(M, g)$ be a compact oriented Riemannian 4-manifold.
It is known that $M$ has nonnegative isotropic curvature if $w_3^\pm \leq s/6$ (see [4]), where $w_3^\pm$ are the largest eigenvalues
of $W^\pm$, respectively. Equation [1.6] and the condition $K_1^\perp \geq \frac{s}{24}$ implies that $ w_1^+ \geq w_1^+ + w_1^-  \geq -s/12$ and so $w_3^+ = -w_1^+ -w_2^+ \leq -2w_1^+ \leq s/6$. Similarly, $w_3^- \leq s/6$ and this proves that $M$ has nonnegative isotropic curvature. Notice that the condition $K_3^\perp \leq \frac{s}{6}$ also implies that $M$ has nonnegative isotropic curvature. On the other hand, it is easy to see that if $M = M_1^2 \times M_2^2$ then $M$ has $K_1^\perp = K_2^\perp = 0$ which contradicts the initial hypotheses. So, this proves that if $M$ is reducible then the universal covering  of $M$ is isometric to product $\mathbb{R} \times N^3$, where $N^3$ is diffeomorphic  to sphere $\mathbb{S}^3$. If $M$ is irreducible then the Theorem 1 follows from principal results in [7] and [1].

Author’s address:
Mathematics Department, Federal University of Bahia,
\\
zipcode: 40170110- Salvador -Bahia-Brazil
\\
Author’s email: ezio@ufba.br or ezioufba@gmail.com
\end{document}